\documentclass[12 pt]{amsart}
\usepackage{amscd,amssymb,amsmath,amsthm}
\input xy
\xyoption{all}
\newtheorem{Theorem}{Theorem}

\newtheorem{Conjecture}{Conjecture}

\newcommand{\proj}{\mathbb{P}}
\newcommand{\Z}{\mathbb{Z}}

\newcommand{\rarr}{\rightarrow}
\newcommand{\oh}{{\mathcal{O}}}
\newcommand{\com}{\mathbb{C}}

\newcommand{\bZ}{{\mathsf{Z}}}
\newcommand{\bF}{{\mathsf{F}}}

\newcommand\udot{^\bullet}

\begin{document}

\title{Maps, sheaves, and $K3$ surfaces}
\author{R. Pandharipande}
\maketitle

\begin{abstract}
The conjectural equivalence of curve counting on Cal\-abi-Yau 3-folds via
stable maps and stable pairs is discussed. By considering
Calabi-Yau 3-folds with $K3$ fibrations, 
the correspondence naturally connects curve and sheaf counting
on $K3$ surfaces. New results and conjectures (with D. Maulik)
about descendent
integration on $K3$ surfaces are announced.
The recent proof of the Yau-Zaslow conjecture
is surveyed . 
\end{abstract}


\baselineskip=18pt

\section{Counting curves}
\subsection{Calabi-Yau 3-folds}
A {\em Calabi-Yau 3-fold} is a nonsingular projective variety{\footnote{All
varieties here are defined over $\com$.} $X$ of
dimension $3$ with trivial first Chern class
$$\wedge^3 T_X \stackrel{\sim}{=} \mathcal{O}_X.$$
Often the triviality of the fundamental group
$$\pi_1(X) = 1$$
is included in the definition. However, for our purposes, $X$ need
not be simply connected.
\subsection{Maps}
Let $C$ be a complete
curve with at worst simple nodes as singularities. 
We do not require $C$ to be connected. The arithmetic genus $g$ of 
$C$ is defined by the Riemann-Roch formula,
$$\chi(C,\mathcal{O}_C)=1-g.$$ 
We view an algebraic map
$$f: C\rightarrow X$$
which is not constant on {\em any} connected component of $C$ 
as 
 parameterizing a subcurve of $X$.
Let
$$\beta= f_*[C] \in H_2(X,\mathbb{Z})$$
be the homology class represented by $f$.
Since $f$ is nonconstant, $\beta \neq 0$.

An automorphism of $f$ is an automorphism of the domain
$$\epsilon:C\rightarrow C$$
satisfying
$f\circ \epsilon = f$.
A map $f$ is {\em stable} \cite{k}
if the automorphism group $\text{Aut}(f)$
is finite.
Infinite automorphisms can come only from contracted
rational and elliptic irreducible components of $C$ incident
to too few nodes.

Let $\overline{M}_g(X,\beta)^\bullet$ denote the moduli space of  
stable maps{\footnote{Usually $\overline{M}_g(X,\beta)$
denotes the moduli space of stable maps with connected
domains. The bullet in our notation indicates 
the possibility of disconnected domains.}}
 from  genus $g$ curves to $X$
representing the class $\beta$.
The moduli space $\overline{M}_g(X,\beta)^\bullet$ is a projective
Deligne-Mumford stack \cite{bm,fp,k}. Certainly, 
$\overline{M}_g(X,\beta)^\bullet$
may be singular, non-reduced, and disconnected.

The most important structure carried by $\overline{M}_g(X,\beta)^\bullet$
is the obstruction theory \cite{Beh,BehFan,LiTian} governing
deformations of maps.
The Zariski tangent space at
$[f]\in \overline{M}_g(X,\beta)^\bullet$ has dimension
$$\text{dim}_\com ( T_{[f]}) = 3g-3 + h^0(C,f^* T_X).$$
The first term on the right corresponds to deformations of the complex
structure of $C$ and the second term to deformations of
 the map with $C$ fixed.{\footnote{We assume, for the given
interpretation
of terms, the domain
 $C$ has no infinitesimal automorphisms.}}
The obstruction space is
$$\text{Obs}_{[f]} = H^1(C,f^*T_X).$$
Formally, we may view the moduli space $\overline{M}_g(X,\beta)^\bullet$
as being cut out by $\text{dim}_\com (\text{Obs}_{[f]})$ 
equations in the tangent space. Hence,  we expect the
dimension of $\overline{M}_g(X,\beta)^\bullet$ to be
\begin{eqnarray*}
\text{dim}_\com^{expected}\Big( \overline{M}_g(X,\beta)^\bullet
 \Big) 
& = & 3g-3 + h^0(C,f^*T_X)-h^1(C,f^*T_X) \\
& = & 
3g-3 + \chi(C, f^*T_X) \\
& = & 3g-3 + \int_C c_1(T_X) + \text{rank}_\com(T_X) (1-g)\\
& = & 0.
\end{eqnarray*}
The third line is by Riemann-Roch. The Calabi-Yau 3-fold
condition is imposed in the fourth line.

Since
all curves in Calabi-Yau 3-folds are expected
to move in 0-dimensional families, we can hope
to count them.
While $\overline{M}_g(X,\beta)^\bullet$ may have large
positive dimensional components, the obstruction theory
provides a virtual class
$$[\overline{M}_g(X,\beta)^\bullet]^{vir} 
\in H_0(\overline{M}_g(X,\beta),\mathbb{Q})$$
in exactly the expected dimension.

Gromov-Witten theory is the curve counting
defined via integration against the virtual
class of $\overline{M}_g(X,\beta)^\bullet$.
The Gromov-Witten invariants of $X$ are
$$N^\bullet_{g,\beta} = \int_{[\overline{M}_g(X,\beta)^\bullet]^{vir}} 1 
\ \in \mathbb{Q}.$$
For fixed nonzero $\beta\in H_2(X,\mathbb{Z})$, let
$$\bZ_{GW,\beta}(u) =\sum_{g} N_{g,\beta}^\bullet \ u^{2g-2} \  \in 
\mathbb{Q}((u)).$$
be the partition function.{\footnote{Sometimes $\bZ_{GW,\beta}(u)$
as defined here
is called the {\em reduced} partition function since 
the constant map contribution are absent.
The constant contributions, calculated in \cite{FP}, will
not arise in our discussion.
}}
Since $\overline{M}_g(X,\beta)^\bullet$ is empty for
$g$ sufficient negative, $\bZ_{GW,\beta}(u)$ is a Laurent
series.

The Gromov-Witten invariants $N_{g,\beta}^\bullet$ should
be viewed as regularized curve counts. 
 The integrals
$N_{g,\beta}^\bullet$ are symplectic invariants. 
A natural idea is to relate the Gromov-Witten invariants
to strict symplecting curve counts after perturbing the
almost complex structure $J$. However, analytic difficulties
arise. A complete understanding of the symplectic 
geometry here has not yet been obtained.

\subsection{Sheaves}
We may also approach curve counting in a Calabi-Yau 3-fold
$X$ via a gauge/sheaf theoretic
approach \cite{pt1,pt2,pt3}.

We would like to construct a moduli space parameterizing
divisors on curves in $X$. If 
the subcurve 
$$\iota:C\subset X$$
is nonsingular, a divisor determines a line bundle $L \rightarrow C$
together with a section $s\in H^0(C,L)$.
The associated torsion sheaf
$$\iota_*(L)= F $$
on $X$ has 1-dimensional support and section $s\in H^0(X,F)$.
However,  for a compact moduli space,
we must allow the support curve $C$ to acquire singularities
and nonreduced structure. The line bundle $L$ must also
be allowed to degenerate.

A {\em pair} $(F,s)$ consists of a sheaf $F$ on $X$  supported
in dimension 1 together with a section $s\in H^0(X,F)$.  A pair $(F,s)$
 is {\em stable} if 
\begin{enumerate}
\item[(i)]
the sheaf $F$ is {pure}, 
\item[(ii)] the section $\oh_X \stackrel{s}{\rightarrow} F$ has 0-dimensional
cokernel.
\end{enumerate}
Purity here simply means
 every nonzero subsheaf of $F$ has support of dimension 1.
As a consequence,
 the scheme theoretic support $C\subset X$ of  $F$ is a 
Cohen-Macaulay curve.
The support of the cokernel (ii) is a finite length subscheme $Z\subset C$.
If the support $C$ is nonsingular, then the stable pair $(F,s)$ is 
uniquely determined by $Z\subset C$.
However, for general $C$, the subscheme $Z$ does not determine $F$ and $s$.

The discrete invariants of a stable pair are the 
holomorphic Euler characteristic
$\chi(F)\in \mathbb{Z}$ and the class{\footnote{$[F]$ 
is the sum of the classes of the irreducible
1-dimensional curves on which $F$ is supported weighted
by the generic length of $F$ on the
curve.}}
 $[F]\in H_2(X,\mathbb{Z})$.
The moduli space $P_n(X,\beta)$ parameterizes stable pairs satisfying
$$\chi(F)=n, \ \ [F]=\beta.$$
After appropriate choices \cite{pt1}, pair stability coincides
with stability arising from geometric
invariant theory \cite{LPPairs1}. The moduli 
space $P_n(X,\beta)$ is a therefore a projective scheme.

To define
invariants, a virtual cycle is required. The usual
deformation theory of pairs is problematic, but
the fixed-determinant deformation theory{\footnote{Every fixed-determinant
deformation of the complex (to any order) is
quasi-isomorphic to a complex arising from a flat deformation
of a stable pair \cite{pt1}. However, the obstruction theory
obtained from derived category deformations differs
from the classical deformation theory of pairs.}}
of the associated \emph{complex}
in the derived category
\begin{equation*}
I\udot=\{\oh_X\stackrel{s}{\rightarrow}F\}\ \in\,D^b(X)
\end{equation*}
is shown in \cite{HT,pt1} to define a perfect obstruction 
theory for $P_n(X,\beta)$ of virtual dimension zero. A virtual cycle 
is then obtained by \cite{Beh,BehFan,LiTian}. The resulting
regularized counts are
$$
P_{n,\beta}=\int_{[P_n(X,\beta)]^{vir}}1 \ \in  \mathbb{Z}.
$$
Let
$$\bZ_{P,\beta}(q) = \sum_{n\in \mathbb{Z}} P_{n,\beta} \ q^n
\ \in \mathbb{Q}((q)).$$
be the partition function. Since $P_n(X,\beta)$ is empty for
$n$ sufficient negative, $\bZ_{P,\beta}(q)$ is a Laurent
series.

The Gromov-Witten invariants $N_{g,\beta}^\bullet$ 
are $\mathbb{Q}$-valued
since $\overline{M}_{g}(X,\beta)^\bullet$ is a Deligne-Mumford
stack, but 
the stable pairs invariants $P_{n,\beta}$
are $\mathbb{Z}$-valued since
$P_n(X,\beta)$ is a scheme.

\section{Correspondence}
\subsection{Two counts}
We have seen there are at least two regularized counting
strategies for curves in Calabi-Yau 3-folds.
While the Gromov-Witten approach may appear closer to
a pure enumerative invariant since
no auxiliary line bundles play a role, in fact the
two theories are equivalent! 

\begin{Conjecture}\label{ooo} For all Calabi-Yau 3-folds $X$ and nonzero curve
classes $\beta \in H_2(X,\mathbb{Z})$,
$$\bZ_{GW,\beta}(u) = \bZ_{P,\beta}(q)$$
after the variable change $-e^{iu}= q$.
\end{Conjecture}

Actually, the variable change $-e^{iu}=q$ is not a priori
well-defined for Laurent series. The issue is addressed by the
following rationality property.

\begin{Conjecture} \label{ttt}
For all Calabi-Yau 3-folds $X$ and 
nonzero curve classes $\beta \in H_2(X,\mathbb{Z})$, the series
$\bZ_{P,\beta}(q)$ is
the Laurent expansion of a rational function 
invariant under $q \leftrightarrow 1/q$.
\end{Conjecture}

In rigid cases,
Conjecture \ref{ooo} implies the contributions of
multiple covers in
Gromov-Witten theory
exactly match the contributions of the divisor choices
on thickened curves
in the theory of stable pairs.
In geometries with moving curves, the meaning
of Conjecture \ref{ooo}
is more subtle.

\subsection{Other counts}
There are other geometric approaches to curve counting
on Calabi-Yau 3-folds. On the map side, a new theory 
of stability has been very recently put forward by 
B. Kim, A. Kresch, and Y.-G. Oh \cite{kko} generalizing the well-known
theory of admissible covers for dimension 1 targets.
On the sheaf side, the older Donaldson-Thomas theory of
ideal sheaf counts \cite{dt,t} is very natural to pursue. 

While Conjectures \ref{ooo} and \ref{ttt} as stated above are from
\cite{pt1}, the relation between Gromov-Witten theory and
sheaf counting was first discovered in the context
of Donaldson-Thomas theory in \cite{mnop1,mnop2}.
Stable pairs appear to be the closest sheaf enumeration
to Gromov-Witten theory. The precise relationship of 
\cite{kko} to the other theories has yet to be discovered, but an
equivalence
almost surely holds.

\subsection{Evidence}
There are three interesting directions which provide evidence
for Conjectures \ref{ooo} and \ref{ttt}.

The first is the study of local Calabi-Yau toric surfaces.{\footnote{A
local Calabi-Yau toric surface is 
the total space of the canonical bundle of any nonsingular, projective,
toric Fano surface.}}
Both the Gromov-Witten and pairs invariants can be
calculated by the virtual localization
formula \cite{GP}. On the Gromov-Witten side, the
topological vertex of \cite{topver,llz2} evaluates the localization
formula. On the stable pairs side, the evaluation is given
by box counting \cite{pt2}. Conjectures 1 and 2 hold.{\footnote{In the
forthcoming paper \cite{moop}, the most general local Calabi-Yau toric
geometry involving the 3-leg vertex is analysed for 
the Gromov-Witten/Donaldson-Thomas correspondence. It is likely
the same path of argument will apply to stable pairs theory
 also.}} However, the
toric examples are necessarily non-compact. The relevance of
toric calculations to compact Calabi-Yau 3-folds is not clear.

The second direction
is progress towards a geometric proof of Conjecture 2.
The obstruction theory for  $P_n(X,\beta)$ 
is self-dual.{\footnote{The obstruction theory is
equipped with a pairing identifying the tangent space
with the dual of the obstruction space \cite{bdt}.}}
By results of K. Behrend \cite{bdt}, there exists a
constructable function
$$\chi^B: P_n(X,\beta) \rightarrow \mathbb{Z}.$$
with integral{\footnote{The integral is defined
by 
$$\int_{P_n(X,\beta)} \chi^B = \sum_{n\in \mathbb{Z}} n \cdot 
\chi_{\text{top}}\Big( (\chi^B)^{-1}(n)\Big)$$
where $\chi_{\text{top}}$ on 
the right is the usual Euler characteristic.}}
 equal to the pairs invariant,
$$\int_{P_n(X,\beta)} \chi^B = P_{n,\beta}.$$
If $P_n(X,\beta)$ is nonsingular, then
$\chi^B=(-1)^{\text{dim}_\com(P_n(X,\beta))}$ is constant
and
$$P_{n,\beta} = (-1)^{\text{dim}_\com(P_n(X,\beta))} 
\chi_{\text{top}}(P_n(X,\beta)).$$
In \cite{pt2}, properties of $\chi^B$ together with an
essential application of Serre duality imply Conjecture 2
for irreducible{\footnote{A class $\beta$
is irreducible if all 1-dimensional subschemes
 representing $\beta$ are
reduced and irreducible.}}
 curves classes $\beta\in H_2(X,\beta)$.
In remarkable recent work of Y. Toda \cite{toda},
using variants of
Bridgeland's stability conditions, wall-crossing
formulas, and Serre duality, the 
rationality of the closely related  series
$$\bZ^\chi_{P,\beta}(q) =\sum_{n\in \mathbb{Z}} 
\chi_{\text{top}}(P_n(X,\beta)) q^n$$
has been proven for {\em all} nonzero classes $\beta\in H_2(X,\beta)$.
A proper inclusion of $\chi^B$ into Toda's argument
should soon lead to a complete proof of Conjecture 2.

The third direction, curve counting on $K3$ surfaces, will
be discussed in Sections \ref{k31} and \ref{k32}. The topic
contains a mix of classical and quantum geometry. While there
has been recent progress, many beautiful open questions
remain.

\section{Curve counting on $K3$ surfaces} \label{k31}
\subsection{Reduced virtual class}
Let $S$ be a $K3$ surface, and let 
$$\beta\in \text{Pic}(S)= H^{1,1}(S,\com) \cap H^2(S,\mathbb{Z})$$
be an nonzero effective curve class.{\footnote{By Poincar\'e duality,
there is a canonical isomorphism $H_2(S,\mathbb{Z}) \stackrel{\sim}{=}
H^2(S,\mathbb{Z})$, so we may view curves classes as
taking values in either theory.}}
By the virtual dimension formula,
\begin{eqnarray*}
\text{dim}_\com^{expected}\Big( \overline{M}_g(S,\beta)^\bullet \Big) 
& = & 
3g-3 + \chi(C, f^*T_S) \\
& = & 3g-3 + 2(1-g)\\
& = & g-1.
\end{eqnarray*}
Let $[f]\in \overline{M}_g(S,\beta)^\bullet$ be a stable map.
There is canonical surjection
\begin{equation}\label{ttrr}
\text{Obs}_{[f]} \rightarrow \com \rightarrow 0
\end{equation}
obtained from the the composition
$$H^1(C,f^*T_S) \stackrel{\sim}{=} H^1(C,f^* \Omega_S)
\stackrel{df}{\rightarrow} H^1(C,\omega_C) \stackrel{\sim}{=}
\com,$$
where the first isomorphism uses 
$$\wedge^2 T_S\stackrel{\sim}{=}\oh_S.$$
The trivial quotient \eqref{ttrr} forces the vanishing
of $[\overline{M}_g(S,\beta)^\bullet]^{vir}$.
However, the obstruction theory can be  modified to
reduce the obstruction space to the kernel of \eqref{ttrr}.
A reduced virtual class
\begin{equation}\label{vt23}
[\overline{M}_g(S,\beta)^\bullet]^{red} \in 
H_{2g}(\overline{M}_g(S,\beta)^\bullet,\mathbb{Q}),
\end{equation}
in dimension 1 greater than expected,
is therefore defined.

By constructing
trivial quotients of $\text{Obs}_{[f]}$
for each connected component of the domain,
the reduced virtual class \eqref{vt23} is easily seen to
be  supported on the locus of curves
with connected domains. Hence, we need only consider
$$\overline{M}_g(S,\beta) \subset \overline{M}_g(S,\beta)^\bullet.$$

We can also consider stable maps from $r$-pointed curves.
The pointed moduli space $\overline{M}_{g,r}(S,\beta)$ 
has a reduced virtual class of dimension $g+r$.

\subsection{Descendents}
\label{gwt} 
The reduced
Gromov-Witten theory of $S$
is defined via integration against $[\overline{M}_{g,r}(S,\beta)]^{red}$.
Let 
$$\text{ev}_i: \overline{M}_{g,r}(S,\beta) \rarr S,$$
$$ L_i \rarr \overline{M}_{g,r}(S,\beta)$$
denote the evaluation maps and cotangent lines bundles associated to
the $r$ marked points.
Let $\gamma_1, \ldots, \gamma_m$ be a basis of $H^*(S,{\mathbb{Q}})$, and
let $$\psi_i = c_1(L_i) \in \overline{M}_{g,r}(S,\beta).$$
The {\em descendent} fields, denoted in the brackets
by $\tau_k(\gamma_j)$, correspond 
to the classes $\psi_i^k \cup \text{ev}_i^*(\gamma_j)$ on the moduli space
of maps. 
Let
\begin{equation}\label{vppq} \Big\langle \tau_{k_1}(\gamma_{l_1}) \cdots
\tau_{k_r}(\gamma_{l_r})\Big\rangle^{red}_{g,\beta} 
= \int_{[\overline{M}_{g,r}(S,\beta)]^{red}} 
\prod_{i=1}^r \psi_i^{k_i} \cup \text{ev}_i^*(\gamma_{l_i})
\end{equation}
denote the descendent
Gromov-Witten invariants. 
Of course \eqref{vppq} vanishes if the integrand does not
match the dimension of the reduced virtual class.

The reduced Gromov-Witten theory is invariant under deformations
of $S$ which preserve $\beta$ as an algebraic class. A standard
argument{\footnote{The
group of isometries of the $K3$ lattice $U^3 \oplus E_8(-1)^2$
acts transitively on elements with fixed norm and divisibility.
The dependence of the reduced Gromov-Witten on only the
norm and divisibility then follows from
the global Torelli Theorem. See \cite{brl}
for a slightly different point of view on the same result.}}
 shows the invariant \eqref{vppq} depends {\em only} on the
norm
$$\langle \beta,\beta \rangle = \int_S \beta\cup \beta$$
and the divisibility of $\beta \in H^2(S,\mathbb{Z})$.

Let us now specialize, for the remainder of Section \ref{tyu}
\label{tyu}, to an elliptically fibered $K3$ surface
$$\nu: S \rightarrow \proj^1$$
with a section. We assume the section and fiber classes
$$\mathbf{s} , \mathbf{f} \in H^2(S,\mathbb{Z})$$
span $\text{Pic}(S)$.
The cone of effective curve classes is
$$V = \{ m\mathbf{s}+  n\mathbf{f} \ | \ m\geq 0,\ n\geq 0,\ 
(m,n)\neq(0,0)
\ \}.$$
Since the norm of $d\mathbf{s}+ dk \mathbf{f}$ is
$2d^2(k-1)$,
effective classes with all divisibilities $d\geq 1$ 
and norms at least $-2d^2$
can be found on $S$. Elementary arguments show the
 integrals \eqref{vppq} vanish in all other
cases.{\footnote{See, for example,
Lemma 2 of \cite{gwnl}.}}

A natural descendent potential function for the
reduced theory of $K3$ surfaces is defined
by
$$\bF^{S}_{g,m}\big(\tau_{k_1}(\gamma_{l_1}) \cdots
\tau_{k_r}(\gamma_{l_r})\big)=
\sum_{n=0}^\infty 
\Big\langle \tau_{k_1}(\gamma_{l_1}) \cdots
\tau_{k_r}(\gamma_{l_r})\Big\rangle^{red}_{g,m\mathbf{s}+ n \mathbf{f} } 
\ q^{m(n-m)}
$$
for $g\geq 0$ and $m\geq 1$.  
The following Conjecture
is made jointly with D. Maulik \cite{mp2}.

\begin{Conjecture} 
$\bF^{S}_{g,m}\big(\tau_{k_1}(\gamma_{l_1})
\cdots
\tau_{k_r}(\gamma_{l_r})\big)
$
is the Fourier expansion in $q$
of a quasi-modular form of level $m^2$ with pole at $q=0$ of
order at most $m^2$.
\end{Conjecture}

By the ring of quasi-modular forms 
of level $m^2$
with possible poles at $q=0$, we mean the algebra generated by 
the Eisenstein series{\footnote{
The Eisenstein series $E_{2k}$ is the modular form
defined by the equation 
$$-\frac{B_{2k}}{4k} E_{2k}(q) = -\frac{B_{2k}}{4k} + \sum_{n\geq 1} 
\sigma_{2k-1}(n) q^{n},$$
where $B_{2n}$ is the $2n^{th}$ Bernoulli number and
$\sigma_n(k)$ is the sum of the $k^{th}$ powers of
the divisors of $n$,
$$\sigma_k(n) = \sum_{i|n} i^k.$$}}
$E_2$ over the ring of modular forms of level $m^2$. 
We have been able to prove Conjecture 3 in the primitive
case $m=1$ by 
relations in the moduli of curves \cite{fpm}, 
degeneration methods \cite{mptop}, 
and  
the elliptic curve results of \cite{op1,op2,op3}. The $m>1$
case appears to require new techniques.

Let $[p]\in H^4(S,\mathbb{Z})$ denote the Poincar\'e dual
of a point.
The simplest of the $K3$ series is the count
of genus $g$ curves passing through $g$ points,
\begin{equation*}
\bF_{g,1}^S ( {\tau_0(p) \cdots \tau_0(p)} ) 
= \eta^{-24}\ \left( -\frac{1}{24}\ 
q\frac{d}{dq} E_2\right)^g 
\end{equation*}
calculated{\footnote{Our indexing conventions
differ slightly from those adopted in \cite{brl}.}} 
by J. Bryan and C. Leung \cite{brl}.
Here
$$\eta(q) = q^{\frac{1}{24}} \prod_{n=1}^\infty (1-q^n)$$
is Dedekind's function.
Similar calculations in genus 1 for $m=2$
have been done in \cite{jleel}.

\subsection{$\lambda_g$ integrals for $K3$ surfaces}\label{bpsk3}
A connection to the enumerative geometry of Calabi-Yau 3-folds
holds for special integrals in the reduced Gromov-Witten
theory of $K3$ surfaces.
Let
\begin{equation}\label{d432}
R_{g,\beta}=\int_{[\overline{M}_g(S,\beta)]^{red}}
(-1)^g \lambda_g
\end{equation}
for effective curve classes $\beta\in H^2(S,\mathbb{Z})$.
Here, the integrand $\lambda_g$ is the top Chern
class of the Hodge bundle 
$$\mathbb{E}_g \rarr \overline{M}_g(S,\beta)$$
with fiber $H^0(C,\omega_C)$ over moduli point
$$[f:C\rarr S]\in \overline{M}_g(S,\beta).$$
See \cite{FP,GP} for a discussion of Hodge classes in
Gromov-Witten theory. 

The integrals \eqref{d432} arise from the following 3-fold
geometry. Let
$$\pi: X \rightarrow \proj^1$$
be a $K3$-fibered Calabi-Yau 3-fold with 
$$\iota :S \stackrel{\sim}{\rightarrow} \pi^{-1}(0) \subset X.$$
 Assume further the family of $K3$ surfaces determined
by $X$ is transverse to the Noether-Lefschetz divisor
in the moduli of $K3$ surface along which  $\beta$
is an algebraic class. Then, the moduli space 
$$\overline{M}_g(S,\beta)\subset
\overline{M}_g(X,\iota_*\beta)$$
is a connected component. The integral \eqref{d432}
is precisely the contribution of $\overline{M}_g(S,\beta)$
to the Gromov-Witten theory of $X$ \cite{gwnl}.
The discussion here may be viewed as an algebraic
analogue of the twistor construction of \cite{brl}.

The definition of the BPS counts{\footnote
{BPS state counts can be extracted from Gromov-Witten
theory via \cite{GV1,GV2}. The counts $r_{g,\beta}$
are conjecturally {\em integers}.}}
 associated to the Hodge integrals
\eqref{d432} is straightforward. Let 
$\alpha\in \text{Pic}(S)$ be a effective primitive class 
The Gromov-Witten potential $F_{{\alpha}}(u,v)$ 
for classes proportional
to ${\alpha}$
is 
$${F}_{{\alpha}}=
\sum_{g=0}^\infty\   \sum_{m=0}^\infty\   R_{g,m\alpha} \ u^{2g-2} 
v^{m{\alpha}}.$$
The
BPS counts $r_{g,m\alpha}$ are uniquely defined 
by the following equation:
\begin{equation*}
F_\alpha =   \ \ \ \sum_{g= 0}^\infty  \ \sum_{m=0}^\infty \
 r_{g,m\alpha} \ u^{2g-2} \sum_{d>0}
\frac{1}{d}\left( \frac{\sin
({du/2})}{u/2}\right)^{2g-2} v^{dm\alpha}. 
\end{equation*}
We have defined BPS counts  
for both primitive and
divisible classes.

The string theoretic calculations of S. Katz, A. Klemm and C. Vafa \cite{kkv}
via heterotic duality yield two conjectures.

\begin{Conjecture} \label{xxx1}
The BPS count $r_{g,\beta}$ depends upon $\beta$ only through the norm
 $\langle \beta,\beta \rangle$.
\end{Conjecture}

Assuming the validity of
Conjecture \ref{xxx1}, let $r_{g,h}$ denote the BPS count associated
to a class $\beta$ satisfying
$$\langle \beta,\beta \rangle = 2h-2.$$
Conjecture \ref{xxx1} is rather surprising from the point
of view of Gromov-Witten theory. The invariants $R_{g,\beta}$
depend a priori 
upon both the norm and the divisibility of $\beta$.

\begin{Conjecture} \label{xxx2}
The BPS counts $r_{g,h}$ are uniquely determined by the
following equation:
$$\sum_{g =0}^\infty
 \sum_{h= 0}^\infty (-1)^g r_{g,h}(\sqrt{z} - \frac{1}
{\sqrt{z}}
)^{2g}q^h =
\prod_{n=1}^\infty \frac{1}{(1-q^n)^{20} (1-zq^n)^2 (1-z^{-1}q^n)^2}.$$
\end{Conjecture}

As a consequence of Conjecture 5, $r_{g,h}$ vanishes if $g>h$ and
$$r_{g,g}=(-1)^g (g+1).$$
The first values are tabulated below:

\vspace{18pt}

\begin{tabular}{|c||ccccc|}
        \hline
\textbf{}
$r_{g,h}$&    $h= 0$ & 1  & 2 & 3 & 4 \\
        \hline \hline
$g=0$ & $1$ & $24$ & $324$ & 
$3200$ &$25650$  \\
1      &  & $-2$ & 
$-54$ & $-800$  & $-8550$      \\
2      & & & $3$ & 
$88$ & $1401$       \\
3      & &  & 
 & $-4$  & $-126$       \\
4      &  &  & 
 &   & 5       \\
       \hline
\end{tabular}

\vspace{18pt}

Conjectures \ref{xxx1} and \ref{xxx2} provide a complete
solution for $\lambda_g$ integrals in the reduced Gromov-Witten
theory of $K3$ surfaces. The answer is compatible with Conjecture 3
as expected since Hodge integrals may be expressed in terms of
descendent integrals \cite{FP}.

\subsection{Stable pairs on $K3$ surfaces}
Let $S$ be a $K3$ surface with
an irreducible class $\beta\in H^2(S,\mathbb{Z})$ satisfying
\begin{equation*}
\langle \beta,\beta \rangle =2h-2,
\end{equation*}
and
let $P_n(S,h)$ denote the associated moduli space of pairs
on $S$.
Consider again the $K3$-fibered Calabi-Yau 3-fold 
$$\pi:X\rightarrow \proj^1.$$
A deformation argument in \cite{pt3} proves
\begin{equation}\label{byy6}
P_n(S,h)\subset P_n(X,\iota_*\beta)
\end{equation}
 is a connected component of the moduli
space of stable pairs of $X$. Moreover, $P_n(S,h)$
is a {\em nonsingular} projective variety \cite{ky,pt3} of
dimension $n+2h-1$.

Let $\Omega_{P}$ be the cotangent bundle of
the moduli space $P_n(S,h)$.
The self-dual 
obstruction theory on $P_n(S,h)$ induced from
the inclusion \eqref{byy6} has obstruction bundle
$\Omega_{P}$. Hence, the contribution of $P_n(S,h)$
to the stable pairs invariants of $X$ is
\begin{eqnarray*}
\bZ^S_{P,h}(y) & =  &
\sum_n \int_{P_n(S,h)} c_{n+2h-1}(\Omega_P)\ y^n \\
& = &
\sum_{n} (-1)^{n+2h-1} e(P_n(S,h))\  y^n.
\end{eqnarray*}
Here, we have written the stable pairs partition
function in the variable $y$ instead of the traditional $q$
since the latter will be reserved for the Fourier
expansions of modular forms.{\footnote{The conflicting uses of
$q$ seem impossible to avoid. The possibilities for confusion are
great.}}

Fortunately, the topological Euler characteristics of 
$P_n(S,h)$ have been calculated by T. Kawai and K. Yoshioka.
By Theorem 5.80 of \cite{ky},
\begin{multline*}
\sum_{h=0}^\infty \sum_{n=1-h}^\infty  
e(P_n(S,h))\  y^n q^h = \\
\left(\sqrt{y}-\frac{1}{\sqrt{y}}\right)^{-2}\
\prod_{n=1}^\infty \frac{1}{(1-q^n)^{20} (1-yq^n)^2(1-y^{-1}q^n)^2} \ .
\end{multline*}
For our pairs invariants, we require the signed Euler characteristics,
$$
\sum_{h=0}^\infty \bZ_h^S(y) \ q^h = 
\sum_{h=0}^\infty \sum_{n=1-h}^\infty (-1)^{n+2h-1} 
e(P_n(S,h))\  y^n q^h. 
$$
Therefore,
$\sum_{h=0}^\infty \bZ_{P,h}^S(y) \ q^h$ equals
$$
-\left(\sqrt{-y}-\frac{1}{\sqrt{-y}}\right)^{-2}\
\prod_{n=1}^\infty \frac{1}{(1-q^n)^{20} (1+yq^n)^2(1+y^{-1}q^n)^2} \ .
$$

\subsection{Correspondence}
We are now in a position to check whether the Katz-Klemm-Vafa
predictions for the $\lambda_g$ integrals in the reduced
Gromov-Witten theory of $S$ are compatible with the above
stable pairs calculations via the maps/pairs correspondence of Conjecture 1.

In the $\beta$ irreducible case, the Gromov-Witten partition
function takes the
form
$$\sum_{h=0}^\infty \bZ_{GW,h}^S(u) \ q^h =
\sum_{g= 0}^\infty \sum_{h=0}^\infty r_{g,h}\ u^{2g-2} 
\left( \frac{\sin
({u/2})}{u/2}\right)^{2g-2}\ q^h.$$
Afte substituting $-e^{iu}=y$, we find
$$\sum_{h=0}^\infty \bZ_{GW,h}^S(y) \ q^h =
\sum_{g= 0}^\infty \sum_{h=0}^\infty (-1)^{g-1} r_{g,h}
 \left(\sqrt{-y}-\frac{1}{\sqrt{-y}}\right)^{2g-2}
\ q^h.$$
By Conjecture \ref{xxx2},
$\sum_{h=0}^\infty \bZ_{GW,h}^S(y) \ q^h$  equals
$$-\left(\sqrt{-y}-\frac{1}{\sqrt{-y}}\right)^{-2}
\prod_{n=1}^\infty \frac{1}{(1-q^n)^{20} (1+yq^n)^2(1+y^{-1}q^n)^2}$$
which is $\sum_{h=0}^\infty \bZ_{P,h}^S(y) \ q^h$.

The maps/pairs correspondence of Conjecture 1 therefore
works perfectly assuming the Katz-Klemm-Vafa prediction
for the reduced Gro\-mov-Witten theory. But
can the Katz-Klemm-Vafa prediction for stable maps be proven?
The answer is yes in genus 0. The proof is our
last topic.

\section{The Yau-Zaslow conjecture}
\label{k32}
\subsection{Genus 0}
The genus 0 parts of Conjectures \ref{xxx1} and \ref{xxx2} 
for $K3$ surfaces were
predicted earlier by S.-T. Yau and E. Zaslow \cite{yauz}.

\vspace{10pt}
\noindent{\bf Conjecture ${\mathbf{4'}}$.}
The BPS count $r_{0,\beta}$ depends upon $\beta$ only through the norm $\langle \beta,\beta \rangle$.
\vspace{10pt}

Let $r_{0,m,h}$ denote the genus 0 BPS count associated
to a class $\beta$  of divisibility $m$ satisfying
$$\langle \beta,\beta \rangle = 2h-2.$$
Assuming Conjecture $4'$ holds, we define
$$r_{0,h}= r_{0,m,h}$$
independent{\footnote{Independence of $m$ holds
when $2m^2$ divides $2h-2$. Otherwise, no such class
$\beta$ exists and $r_{0,m,h}$ is defined
to vanish.}} of $m$.

\vspace{10pt}
\noindent{\bf Conjecture ${\mathbf{5'}}$.}
The BPS counts $r_{0,h}$ are uniquely determined by
\begin{equation} \label{gvre}
 \sum_{h\geq 0} r_{0,h}\ q^{h-1} = q^{-1}{\prod_{n=1}^\infty (1 - q^n)^{-24}}.
\end{equation}
\vspace{10pt}

A mathematical derivation of the Yau-Zaslow conjectures
for primitive classes $\beta$ 
via 
Euler characteristics of compactified Jacobians 
following \cite{yauz} 
can be found in \cite{beu,xc,fgd}.
The Yau-Zaslow formula \eqref{gvre}
 was proven via Gromov-Witten theory for
primitive classes $\beta$ by J. Bryan and C. Leung \cite{brl}. 
An early calculation by A. Gathmann \cite{Gat} for a class
$\beta$ of
divisibility 2 was important
for the correct formulation of the conjectures.  
Conjectures $4'$ and $5'$ have been proven in the divisibility 2 case
by J. Lee and C. Leung \cite{ll} and B. Wu \cite{wwuu}.

The main result of the paper \cite{kmps} with A. Klemm, D. Maulik,
and E. Scheidegger  
is a proof of Conjectures ${4'}$ and ${5'}$
in all cases.

\begin{Theorem} \label{yzz}
The Yau-Zaslow conjectures hold for all nonzero effective classes 
$\beta\in \text{\em Pic}(S)$ on a $K3$ surface $S$.
\end{Theorem}

The proof, using the connection to 
Noether-Lefschetz theory \cite{gwnl}, mirror symmetry, and
modular form identities, is surveyed in Sections \ref{nnll} -\ref{bb23}.

\subsection{Noether-Lefschetz theory} \label{nnll}
\subsubsection{$K3$ lattice}
Let $S$ be a $K3$ surface.
 The second cohomology of $S$ is a rank 22 lattice
with intersection form 
\begin{equation}\label{ccet}
H^2(S,\mathbb{Z}) \stackrel{\sim}{=} U\oplus U \oplus U \oplus E_8(-1) \oplus E_8(-1)
\end{equation}
where
$$U
= \left( \begin{array}{cc}
0 & 1 \\
1 & 0 \end{array} \right)$$
and 
$$E_8(-1)=  \left( \begin{array}{cccccccc}
 -2&    0 &  1 &   0 &   0 &   0 &   0 & 0\\
    0 &   -2 &   0 &  1 &   0 &   0 &   0 & 0\\
     1 &   0 &   -2 &  1 &   0 &   0 & 0 &  0\\
      0  & 1 &  1 &   -2 &  1 &   0 & 0 & 0\\
      0 &   0 &   0 &  1 &   -2 &  1 & 0&  0\\
      0 &   0&    0 &   0 &  1 &  -2 &  1 & 0\\ 
      0 &   0&    0 &   0 &   0 &  1 &  -2 & 1\\
      0 & 0  & 0 &  0 & 0 & 0 & 1& -2\end{array}\right)$$
is the (negative) Cartan matrix. The intersection form \eqref{ccet}
is even.

\subsubsection{Lattice polarization}
A primitive{\footnote{A class in $H^2(S,\mathbb{Z})$ of divisibility
1 is
{\em primitive}.}} class 
$L\in \text{Pic}(S)$ is a {\em quasi-polarization}
if
$$\langle L,L \rangle >0  \ \ \ \text{and} \ \  \ \langle L,[C]\rangle
 \geq 0 $$
for every curve $C\subset S$.
A sufficiently high tensor power $L^n$
of a quasi-polarization is base point free and determines
a birational morphism
$$S\rightarrow \tilde{S}$$
contracting A-D-E configurations of $(-2)$-curves on $S$.
Hence, every quasi-polarized $K3$ surface is algebraic.

Let $\Lambda$ be a fixed rank $r$  
primitive{\footnote{An embedding
of lattices is primitive if the quotient is torsion free.}}
embedding
\begin{equation*} 
\Lambda \subset U\oplus U \oplus U \oplus E_8(-1) \oplus E_8(-1)
\end{equation*}
with signature $(1,r-1)$, and
let 
$v_1,\ldots, v_r \in \Lambda$ be an integral basis.
The discriminant is
$$\Delta(\Lambda) = (-1)^{r-1} \det
\begin{pmatrix}
\langle v_{1},v_{1}\rangle & \cdots & \langle v_{1},v_{r}\rangle  \\
\vdots & \ddots & \vdots \\
\langle v_{r},v_{1}\rangle & \cdots & \langle v_{r},v_{r}\rangle 
\end{pmatrix}\ .$$
The sign is chosen so $\Delta(\Lambda)>0$.

A {\em $\Lambda$-polarization} of a $K3$ surface $S$   
is a primitive embedding
$$j: \Lambda \rightarrow \mathrm{Pic}(S)$$ 
satisfying two properties:
\begin{enumerate}
\item[(i)] the lattice pairs 
$\Lambda \subset U^3\oplus E_8(-1)^2$ and
$\Lambda\subset 
H^2(S,\mathbb{Z})$ are isomorphic
 via an isometry which restricts to the identity on $\Lambda$,
 \item[(ii)]
$\text{Im}(j)$ contains
a {quasi-polarization}. 
\end{enumerate}
By (ii), every $\Lambda$-polarized $K3$ surface is algebraic.

The period domain $M$ of Hodge structures of type $(1,20,1)$ on the lattice 
$U^3 \oplus E_8(-1)^2$   is
an analytic open set
 of the 20-dimensional  nonsingular isotropic 
quadric $Q$,
$$M\subset Q\subset \proj\big(    (U^3 \oplus E_8(-1)^2 )    
\otimes_\Z \com\big).$$
Let $M_\Lambda\subset M$ be the locus of vectors orthogonal to 
the entire sublattice $\Lambda \subset U^3 \oplus E_8(-1)^2$.

Let $\Gamma$ be the isometry group of the lattice 
$U^3 \oplus E_8(-1)^2$, and let
 $$\Gamma_\Lambda \subset \Gamma$$ be the
subgroup  restricting to the identity on $\Lambda$.
By global Torelli,
the moduli space $\mathcal{M}_{\Lambda}$ 
of $\Lambda$-polarized $K3$ surfaces 
is the quotient
$$\mathcal{M}_\Lambda = M_\Lambda/\Gamma_\Lambda.$$
We refer the reader to \cite{dolga} for a detailed
discussion.

\subsubsection{Families}

Let $X$ be a 
compact 3-dimensional complex manifold equipped with
holomorphic line bundles
$$L_1, \ldots, L_r  \ \rightarrow X$$
 and a holomorphic map
 $$\pi: X \rarr C$$
to a nonsingular complete curve.

The tuple $(X,L_1,\ldots, L_r, \pi)$ is a
{\em 1-parameter family of nonsingular $\Lambda$-polarized
$K3$ surfaces}
if 
\begin{enumerate}
\item[(i)] the fibers $(X_\xi, L_{1,\xi}, \ldots, L_{r,\xi})$
are $\Lambda$-polarized $K3$ surfaces via
$$v_i \mapsto L_{i,\xi}$$
for every $\xi\in C$,
\item[(ii)] there exists a $\lambda^\pi\in \Lambda$
which is a quasi-polarization of all fibers of $\pi$ 
simultaneously.
\end{enumerate}
The family $\pi$ yields a morphism,
$$\iota_\pi: C \rarr \mathcal{M}_{\Lambda},$$
to the moduli space of $\Lambda$-polarized $K3$ surfaces. 

Let $\lambda^{\pi}= \lambda^\pi_1 v_1+\dots +\lambda^\pi_r v_r$.
A vector $(d_1,\ldots,d_r)$ of integers is {\em positive} if
$$\sum_{i=1}^r \lambda^\pi_i d_i >0.$$
If $\beta \in \text{Pic}(X_\xi)$ has intersection numbers
$$d_i = \langle L_{i,\xi},\beta \rangle,$$
then $\beta$ has positive degree with respect to the
quasi-polarization if and only if  $(d_1,\dots,d_r)$
is positive.

\subsubsection{Noether-Lefschetz divisors}
Noether-Lefschetz numbers are defined in \cite{gwnl}
by the intersection of $\iota_\pi(C)$ with Noether-Lefschetz 
divisors in $\mathcal{M}_\Lambda$.
We briefly review the definition of the 
Noether-Lefschetz divisors.

Let $(\mathbb{L}, \iota)$ be a rank $r+1$
lattice  $\mathbb{L}$
with an even symmetric bilinear form $\langle,\rangle$ and a primitive embedding
$$\iota: \Lambda \rightarrow \mathbb{L}.$$
Two data sets 
$(\mathbb{L},\iota)$ and $(\mathbb{L}',  \iota')$
are isomorphic if there is an isometry which restricts to identity on $\Lambda$.
The first invariant of the data $(\mathbb{L}, \iota)$ is
the discriminant
 $\Delta \in \mathbb{Z}$ of 
$\mathbb{L}$.

An additional invariant of $(\mathbb{L}, \iota)$ can be 
obtained by considering 
any vector $v\in \mathbb{L}$ for which
\begin{equation}\label{ccff} 
\mathbb{L} = \iota(\Lambda) \oplus \mathbb{Z}v.
\end{equation}
The pairing
$$\langle v, \cdot \rangle : \Lambda \rightarrow \mathbb{Z}$$
determines an element of $\delta_v\in \Lambda^*$.
Let 
$G = \Lambda^{*}/\Lambda$
be quotient defined via the injection
$\Lambda \rightarrow \Lambda^*$
 obtained from the pairing $\langle,\rangle$ on $\Lambda$.
The group $G$ is abelian of order equal to the 
discriminant $\Delta(\Lambda)$.
The image 
$$\delta \in G/\pm$$
of $\delta_v$ is easily seen to be independent of $v$ satisfying 
\eqref{ccff}. The invariant $\delta$ is the {\em coset} of $(\mathbb{L},\iota)$

By elementary arguments, two data sets $(\mathbb{L},\iota)$ and $(\mathbb{L}',\iota')$
of  rank $r+1$ are isomorphic if and only if the discriminants and cosets are
equal.

Let $v_1,\ldots, v_r$ be an integral basis of $\Lambda$ as before.
The pairing of $\mathbb{L}$ 
with respect to an extended basis $v_{1}, \dots, v_{r},v$
is encoded in the matrix
$$\mathbb{L}_{h,d_{1},\dots,d_{r}} = 
\begin{pmatrix}
\langle v_{1},v_{1}\rangle & \cdots & \langle v_{1},v_{r}\rangle & d_{1} \\
\vdots & \ddots & \vdots & \vdots\\
\langle v_{r},v_{1}\rangle & \cdots & \langle v_{r},v_{r}\rangle & d_{r}\\
d_{1} & \cdots & d_{r} & 2h-2
\end{pmatrix}.$$
The discriminant is
$$\Delta(h,d_{1},\dots,d_{r}) 
= (-1)^r\mathrm{det}(\mathbb{L}_{h,d_{1},\dots,d_{r}}).$$
The coset $\delta(h, d_{1},\dots,d_{r})$ is represented by the functional
$$v_i \mapsto d_i.$$

The Noether-Lefschetz divisor $P_{\Delta,\delta} \subset \mathcal{M}_{\Lambda}$
is the closure of the locus of $\Lambda$-polarized $K3$ surfaces $S$ for which
$(\mathrm{Pic}(S),j)$ has rank $r+1$, discriminant $\Delta$, and coset $\delta$.
By the Hodge index theorem, $P_{\Delta,\delta}$ is empty unless $\Delta > 0.$

Let $h, d_{1}, \dots, d_{r}$ determine a positive discriminant
$$\Delta(h,d_{1},\dots,d_{r}) > 0.$$  The Noether-Lefschetz divisor
$D_{h, (d_{1},\dots,d_{r})}\subset \mathcal{M}_{\Lambda}$ is defined by 
the weighted sum
$$D_{h,(d_{1},\dots,d_{r})} 
= \sum_{\Delta,\delta} m(h,d_1,\dots,d_r|\Delta,\delta)\cdot[P_{\Delta,\delta}]$$
where the multiplicity $m(h,d_1,\dots,d_r|\Delta,\delta)$ is the number
of elements $\beta$ of the lattice $(\mathbb{L},\iota)$ of 
type $(\Delta,\delta)$ satisfying
\begin{equation}\label{34f}
\langle \beta, \beta \rangle = 2h-2,\ \  \langle \beta, v_{i}\rangle = d_{i}.
\end{equation}
If the multiplicity is nonzero, then $\Delta | \Delta(h,d_{1},\dots,d_{r})$ so only 
finitely many divisors appear in the 
above sum.

If $\Delta(h,d_{1},\dots,d_{r}) = 0$, the divisor $D_{h,(d_1,\dots,d_r)}$
has an alternate definition.
The tautological
line bundle $\mathcal{O}(-1)$ is $\Gamma$-equivariant on the period domain
$M_\Lambda$ and descends to the {\em Hodge line bundle} 
$$\mathcal{K} \rightarrow \mathcal{M}_{\Lambda}.$$
We define
$D_{h,(d_{1},\dots,d_{r})} = \mathcal{K}^{*}$.
See \cite{gwnl} for an alternate view of degenerate intersection.

If $\Delta(h,d_{1},\dots,d_{r}) < 0$, the divisor 
$D_{h,(d_1,\dots,d_r)}$ on $\mathcal{M}_{\Lambda}$ is defined to vanish
by the Hodge index theorem.

\subsubsection{Noether-Lefschetz numbers}
Let $\Lambda$ be a lattice of discriminant $l=\Delta(\Lambda)$, and
let $(X,L_1,\ldots,L_r,\pi)$ be 
a 1-parameter family of $\Lambda$-polarized $K3$ surfaces.
The Noether-Lefschetz number $NL^\pi_{h,d_1,\dots,d_r}$ is
the  classical intersection
product
\begin{equation}\label{def11}
NL^\pi_{h,(d_1,\dots,d_r)} =\int_C \iota_\pi^*[D_{h,(d_1,\dots,d_r)}].
\end{equation}

Let $\mathrm{Mp}_{2}(\mathbb{Z})$ be
the metaplectic double cover of $SL_{2}(\mathbb{Z})$. 
There is a canonical representation \cite{borch}
associated to $\Lambda$,
$$\rho_{\Lambda}^{*}: 
\mathrm{Mp}_{2}(\mathbb{Z}) \rightarrow \mathrm{End}(\mathbb{C}[G]).$$
The full set of Noether-Lefschetz numbers
$NL^\pi_{h,d_1,\dots,d_r}$ defines a vector valued
modular form
 $$\Phi^{\pi}(q) = \sum_{\gamma\in G} \Phi^{\pi}_{\gamma}(q)v_{\gamma} \in \com[[q^{\frac{1}{2l}}]]
\otimes \com[G],$$
of weight $\frac{22-r}{2}$ and type $\rho_\Lambda^*$
by results{\footnote{While the results of the papers \cite{borch, kudmil}
have considerable overlap, we will follow the point of view of Borcherds.}}
 of Borcherds and Kudla-Millson \cite{borch,kudmil}.
The Noether-Lefschetz numbers are the coefficients{\footnote{If $f$ is a series in $q$,
 $f[k]$ denotes the coefficient of $q^k$.}}
 of the components of 
$\Phi^\pi$,
$$NL^{\pi}_{h,(d_1,\dots,d_r)} = \Phi^{\pi}_{\gamma}\left[ \frac{\Delta(h,d_1,\dots,d_r)}{2l}\right]$$
where $\delta(h,d_1,\dots,d_r) = \pm\gamma$.
The modular form results significantly constrain the Noether-Lefschetz numbers.

\subsubsection{Refinements}
If  $d_1,\ldots,d_r$ do not simultaneously
vanish, refined Noether-Lefschetz divisors
are defined.
If $\Delta(h,d_1,\dots,d_r)>0$, 
$$D_{m,h,(d_1,\dots,d_r)}
\subset D_{h,(d_1,\dots,d_r)}$$ is defined
by requiring the class $\beta \in \text{Pic}(S)$ to satisfy \eqref{34f} and
have divisibility $m>0$. If $\Delta(h,d_1,\dots,d_r)=0$, then
$$D_{m,h,(d_1,\dots,d_r)}=D_{h,(d_1,\dots,d_r)}$$
if $m>0$ is the greatest common divisor of $d_1,\ldots,d_r$
and 0 otherwise.

Refined
Noether-Lefschetz numbers are defined by
\begin{equation}\label{def112}
NL^\pi_{m,h,(d_1,\dots,d_r)} =\int_C \iota_\pi^*[D_{m,h,(d_1,\dots,d_r)}].
\end{equation}
The full 
set of Noether-Lefschetz numbers $NL^\pi_{h,(d_1,\dots,d_r)}$ is
easily shown in \cite{kmps}
to determine the refined numbers $NL^\pi_{m,h,(d_1,\dots,d_r)}$.

\subsection{Three theories}
The main geometric idea in the proof of Theorem 1  
is the relationship  of three theories associated to   
a 1-parameter family $$\pi:X \rightarrow C$$
 of $\Lambda$-polarized $K3$ surfaces:
\begin{enumerate}
\item[(i)] the Noether-Lefschetz numbers  of $\pi$,
\item[(ii)] the genus 0 Gromov-Witten invariants of $X$,
\item[(iii)] the genus 0 reduced Gromov-Witten invariants of the
$K3$ fibers. 
\end{enumerate}
The Noether-Lefschetz numbers (i) are classical intersection
products while the Gromov-Witten invariants (ii)-(iii) 
are quantum in origin. 
For (ii),
we view the theory in terms  the
Gopakumar-Vafa invariants
\cite{GV1,GV2}.

Let $n_{0,(d_1,\dots,d_r)}^X$ denote the Gopakumar-Vafa invariant of $X$
in genus $0$ for $\pi$-vertical curve classes of degrees $d_1,\dots,d_r$
with respect to the line bundles $L_1,\dots, L_r$. Let
$r_{0,m,h}$ denote the reduced $K3$ invariant.
The following result is proven{\footnote{The result of
the \cite{gwnl} is stated in the rank $r=1$ case, but the
argument is identical for arbitrary $r$.}} in \cite{gwnl}
by a comparison of
the reduced and usual deformation theories of maps of curves
to
the $K3$ fibers of $\pi$.

\begin{Theorem} \label{ffc}
 For degrees $(d_1,\dots,d_r )$ positive with respect to the
quasi-polariza\-tion $\lambda^\pi$,
$$n_{0,(d_1,\dots,d_r)}^X= \sum_{h=0}^\infty \sum_{m=1}^{\infty}
r_{0,m,h}\cdot  NL_{m,h,(d_1,\dots,d_r)}^\pi.$$
\end{Theorem}

\subsection{The STU model}
The STU model{\footnote{The model has been studied in physics
since the
80's. The
letter $S$ stands for the dilaton and $T$ and $U$
label the torus moduli in the heterotic string. The STU model
was an important example for the
duality between type IIA and heterotic strings
formulated in \cite{kv} and has been intensively studied
 \cite{hm1,hm2,germans,klm,MM}.}}
is a particular nonsingular projective
Calabi-Yau 3-fold $X$
equipped with a fibration 
\begin{equation}\label{fibbw}
\pi: X \rightarrow \proj^1.
\end{equation}
Except for 528 points $\xi \in \proj^1$,
the fibers 
$$X_\xi= \pi^{-1}(\xi)$$ are nonsingular
elliptically fibered $K3$ surfaces.
The 528 singular fibers $X_\xi$ have exactly 1 ordinary
double point singularity each.

The 3-fold $X$ is constructed as a nonsingular
anticanonical section of 
the nonsingular projective toric 4-fold $Y$ defined by
10 rays with primitives
\vspace{5pt}
\begin{center}
\begin{tabular}{lll}
$\rho_1= (1,0,2,3)\ \ \ $ & $\rho_2=(-1,0,2,3)\ \ \ $  & \\
$\rho_3=(0,1,2,3)$ &  $\rho_4=(0,-1,2,3)$  & \\
$\rho_5= (0,0,2,3)$ &  $\rho_6=(0,0,-1,0)$ &  $\rho_7=(0,0,0,-1)$\\
$\rho_8 = (0,0,1,2)$  & $\rho_9=(0,0,0,1)$ &  $\rho_{10}=(0,0,1,1)$.
\end{tabular}
\end{center}
\vspace{15pt}
The Picard rank of $Y$ is 6. The fibration \eqref{fibbw}
is obtained from a nonsingular toric fibration
$$\pi^Y: Y \rightarrow \proj^1.$$
The image of
$$\text{Pic}(Y)\rightarrow \text{Pic}(X_\xi)$$
determines a rank 2 sublattice of each fiber $\text{Pic}(X_\xi)$
with intersection form
\begin{equation*}
 \Lambda= \left( \begin{array}{cc}
0 & 1 \\
1 & 0 \end{array} \right)\ .
\end{equation*}
Let $L_1, L_2 \rightarrow X$ denote line bundles
which span the standard basis of the form $\Lambda$ after restriction.

Strictly speaking, the tuple $(X,L_1,L_2,\pi)$ is
not a $1$-parameter family of $\Lambda$-polarized
$K3$ surfaces.
The only failing
is the 528 singular fibers of $\pi$.
Let
$$\epsilon:C \stackrel{2-1}{\longrightarrow} \proj^1$$
be a hyperelliptic curve branched over the 528 points of
$\proj^1$ corresponding to the singular fibers of $\pi$.
The family
$$\epsilon^*(X) \rightarrow C$$
has 3-fold double point singularities over the 528 nodes
of the fibers of the original family.
Let
$$\widetilde{\pi}: \widetilde{X} \rightarrow C$$
be obtained from a small resolution
$$\widetilde{X} \rightarrow \epsilon^*(X).$$
Let $\widetilde{L}_i\rightarrow \widetilde{X}$ be the
pull-back of $L_i$ by $\epsilon$.
The data $$(\widetilde{X},  \widetilde{L}_1,
\widetilde{L}_2,\widetilde{\pi})$$ determine a 1-parameter
family of $\Lambda$-polarized
$K3$ surfaces, see Section 5.3 of \cite{gwnl}.
The simultaneous quasi-polarization is obtained from
the projectivity of $X$.

\subsection{Proof of Theorem \ref{yzz}}
Theorem \ref{yzz} is proven in \cite{kmps} by studying 
Theorem 2 applied to the STU model. There 
are four basic steps:
\begin{enumerate}
\item[(i)] 
The modular form \cite{borch,kudmil} determining the intersections
of the base curve 
with the Noether-Lefschetz divisors is calculated.
For the STU model, the modular form has vector dimension 1 and is
proportional to the product $E_4 E_6$ of Eisenstein series.
\item[(ii)]
Theorem \ref{ffc} is used to show the 3-fold BPS counts
$n_{0,(d_1,d_2)}^{\widetilde{X}}$ then {\em determine} all the reduced $K3$ 
invariants $r_{0,m,h}$. Strong use is made of the rank 2 lattice of
the STU model. 
\item[(iii)]
The BPS counts $n_{0,(d_1,d_2)}^{\widetilde{X}}$ are calculated
via mirror symmetry. Since the STU model is realized as a Calabi-Yau
complete intersection in a nonsingular toric variety, the genus 0
Gromov-Witten invariants are obtained after proven mirror transformations
from hypergeometric series \cite{giv1,giv2,lly}. The Klemm-Lerche-Mayr
identity, proven in \cite{kmps}, shows the
invariants $n_{0,(d_1,d_2)}^{\widetilde{X}}$ 
are themselves related to modular forms.
\item[(iv)] 
Theorem \ref{yzz}  then follows from the Harvey-Moore
identity which simultaneously relates the modular structures of
$$n_{0,(d_1,d_2)}^{\widetilde{X}}, \ \ r_{0,m,h}, \ \ \text{and}\ \  
NL^{\widetilde{\pi}}_{m,h,(d_1,d_2)}$$
 in the
form specified by Theorem \ref{ffc}.
\end{enumerate}

The Harvey-Moore identity of part (iv) is simple to state.
Let
$$f(\tau)= \frac{E_4(\tau)E_6(\tau)}{\eta(\tau)^{24}} =
\sum_{n=-1}^\infty c(n) q^n $$
where $q=\exp(2\pi i \tau)$.
Then, 
\begin{equation}
\label{harmoo}
\frac{f(\tau_1) E_4(\tau_2)}{j(\tau_1)-j(\tau_2)}
= \frac{q_1}{q_1-q_2} + E_4(\tau_2) - \sum_{d,k,\ell>0}
\ell^3 c(k\ell)\ q_1^{kd}q_2^{\ell d}.
\end{equation}
Equation \eqref{harmoo} was 
conjectured in \cite{hm1} and proven by D. Zagier --- the proof
is presented
in Section 4 of \cite{kmps}.

The strategy of the proof of 
the Yau-Zaslow conjectures
is special to genus 0. Much less is known in higher
genus.  
For genus $1$,  
the Katz-Klemm-Vafa conjectures follow for all classes on $K3$ surfaces 
from the Yau-Zaslow conjectures
via 
the boundary relation for $\lambda_1$ in the moduli of
elliptic curves. 
In genus 2 and 3, A. Pixton \cite{pix} has proven the Katz-Klemm-Vafa
formula for primitive classes using boundary relations for
$\lambda_2$ and $\lambda_3$ on $\overline{M}_2$ and
$\overline{M}_3$ respectively.
New ideas will be required for
a complete proof of Conjectures 4 and 5.

\label{bb23}

\section{Acknowledgments}
The paper accompanies my lecture at the Clay research conference
in Cambridge, MA 
in May 2008. The discussion of the enumerative geometry
of stable pairs in Section 1 and 2 reflects joint work with 
R. Thomas. The  study of descendent integrals in 
the reduced Gromov-Witten theory of
$K3$ surfaces in Section 3 is joint work with D. Maulik. The proof
of the Yau-Zaslow conjectures reported in Section 4 is 
joint work with A. Klemm, D. Maulik, and E. Scheidegger.
My research is partially supported by 
NSF grant DMS-0500187.

Many of the ideas discussed here are valid in more general contexts.
For example, the stable maps/pairs correspondence is conjectured in
\cite{pt1} for all 3-folds --- the Calabi-Yau condition is not
necessary. The $K3$ study can be pursued along similar lines for
abelian surfaces, see \cite{ble} for a start. 
The Enriques surface is a close cousin \cite{mp}.

\end{document}